\documentclass[11pt]{article}
\usepackage{amssymb}
\usepackage{enumerate}
\usepackage{enumitem}
\usepackage{mathrsfs}
\usepackage{latexsym,amsmath,amssymb,amsfonts,epsfig,graphicx,cite,psfrag}
\usepackage{eepic,color,colordvi,amscd}
\usepackage{verbatim}
\usepackage{appendix}

\usepackage{float}
\usepackage{graphicx}
\usepackage{cite}
\usepackage{caption}
\usepackage{subcaption}

\newcommand{\qed}{\hfill $\Box $}
\newcommand{\pf}{\noindent {\bf Proof.} }
\newtheorem{theorem}{Theorem}
\newtheorem{lemma}[theorem]{Lemma}

\newcommand{\ex}{{\rm ex}}
\newcommand{\G}{{\cal G}}

\begin{document}
	
	\title{Anti-Ramsey Number of Friendship Graphs\thanks{Supported by National Key R\&D Program of China under grant No.2023YFA1010203 and  the National Natural
Science Foundation of China under grant No.12271425}}

\author{
 Wenke Liu, Hongliang Lu\footnote{Corresponding author: luhongliang215@sina.com}  and Xinyue Luo\\School of Mathematics and Statistics\\
Xi'an Jiaotong University,
Xi'an, Shaanxi 710049, China\\
\smallskip\\
}
\date{}


\date{}

\maketitle
	
%
%
%
	
	
	\date{}
	
	\maketitle
	
	\begin{abstract}
		An edge-colored graph is called \textit{rainbow graph} if all the colors on its edges are distinct. For a given positive integer $n$ and a family of graphs $\mathcal{G}$, the anti-Ramsey number $ar(n, \mathcal{G})$ is the smallest number of colors $r$ required to ensure that, no matter how the edges of the complete graph $K_n$ are colored using exactly $r$ colors, there will always be a rainbow copy of some graph $G$ from the family $\mathcal{G}$.
A friendship graph $F_k$ is the graph obtained  by combining $k$ triangles that share a common vertex. In this paper, we determine the anti-Ramsey number $ar(n, \{F_k\})$ for large values of $n$. Additionally, we also determine the $ar(n, \{K_{1,k}, kK_2\}$, where $K_{1,k}$ is a star graph with $ k$ vertices and $kK_2$ is a matching of size $k$.
		
	\end{abstract}
	\begin{flushleft}
	\end{flushleft}

	\section{Introduction}
	
	We consider finite graphs without multiple edges or loops. Let  $G$ be a graph with vertex set $V (G)$ and edge set $E(G)$.  The number of edges in $G$ is denoted by $e(G)$, i.e., $e(G):=|E(G)|$. We use $|G|$ to denote the number of vertices of $G$.
	For $X\subseteq V (G)$,  we denote the graph obtained by removing all vertices in $X$ and their incident edges from $G$ as $G-X$. When $X$ consists of a single vertex $x$, we simplify the notation to $G-x$ instead of $G-\{x\}$. For any edge set $Y$, we define $G-Y$ as the graph obtained by removing all edges in $Y$ from $G$ and $V(Y)$ is defined as the set of all vertices that are incident to at least one edge in $Y$.
 $G + Y$ is defined as  graph with vertex set $V(G)\cup V(Y)$ and the edge set $E(G)\cup Y$.
   When $Y$ is a single edge $e$, we also use the notation $G+e$ to mean $G +\{e\}$. A subgraph $H$ of $G$ is called an \emph{induced subgraph} if every pair of vertices in $H$ that are adjacent in $G$ remain adjacent in $H$.  An \emph{edge-induced subgraph} is a subset of the edges of $G$ together with any vertices that are their endpoints.
If $H$ and $G$ are graphs and $G$ does not contain $H$ as an  subgraph, we shall say that $G$ is \emph{$H$-free}. Given integers $ n\ge p\ge 1 $, let $ T_{n,p} $ denote the \textit{Tur\'an graph}, i.e., the complete $ p $-partite graph on $ n  $ vertices where each part has either $ \lfloor n/p\rfloor $ or $ \lceil n/p \rceil  $ vertices and the edge set consists of all pairs joining distinct parts. Let $t_p(n)$ denote the number of edges in $ T_{n,p}$. Let $K_n$ denote the complete graph with $ n $ vertices. The  set $\{1,2, \cdots ,n \}$ is denoted by $[n]$.
	
	For a vertex $ x\in V(G) $, the neighborhood  of $ x $ in $G$ is denoted by $N_G(x)$. The \textit{degree} of $ x $ in $ G $, denoted by $ d_G(x )$, is the size of $ N_G(x) $. We use $\delta(G)$ and $ \Delta(G) $ to denote the \textit{minimum} and \textit{maximum degrees} of the vertices of $ G $. A \textit{$k$-regular graph} is a graph in which each vertex has degree $k$. If all vertices have degree $k$ but exactly one vertex has degree $k-1$ in a graph $G$, we call $G$ a \textit{nearly $k$-regular graph}.
	
	For non-empty subset $ X,Y \subseteq V(G) $, let $G[X]$ denote the subgraph induced by $X$ and let $E_G (X,Y)$ denote the set of all edges in $G$, with one endpoint in $X$ and one endpoint in $Y$. 
	A \textit{matching} in $G$ is a set of edges from $ E(G) $, no two of which share a common vertex, and
a \textit{perfect matching} in $ G $ is a matching in which every vertex of $G$ is incident to exactly one edge of the matching. The \textit{matching number} of $G$, denoted by $ \nu(G) $, is the maximum number of edges in a matching in $G$.  Let $G(X,Y)$ denote a bipartite graph with bipartition $X $ and $ Y$. A \textit{near perfect matching} in a graph $G$ is a matching of $G$ covering all but exactly one vertex.
	A graph $G$ is said to be \textit{factor-critical} if $G - v$ has a perfect matching for every $v \in V(G)$.
	
	For two vertex-disjoint graphs $G,H$, the \textit{union} of $G$ and $H$ is the graph $G \cup H$ with vertex set $V(G) \cup V(H)$ and edge set $E(G) \cup E(H)$. The \textit{join} of $G$ and $H$, denoted by $G \lor H$, is the graph obtained from $G \cup H$ by adding edges joining every vertex of $G$ to every vertex of $H$.
	
	An $r$-\textit{edge-coloring} of a graph is an assignment of $r$ colors to the edges of the graph. An \textit{exactly $r$-edge-coloring} of a graph is an $r$-edge-coloring that uses all $r$ colors. An edge-colored graph is called \textit{rainbow} if all edges have distinct colors. Given a positive integer $ n $ and a graph set $ \mathcal{G} $ , the \textit{anti-Ramsey number} $ ar(n,\mathcal{G}) $ is the minimum number of colors $r$ such that each edge-coloring of $K_n$  using exactly
$r$ colors will include a rainbow copy of some graph
$G$ from the set $\mathcal{G}$.
	When $\mathcal{G}=\{G\}$, we denote $ar(n,\mathcal{G})$  by $ ar(n,G) $ for simplicity. The \textit{Tur\'an number} $ \ex(n,G) $ is the maximum number of edges of a graph on $ n $ vertices containing no subgraph isomorphic to $ G $. 	 A
	graph $H$ on $n$ vertices with $\ex(n, G)$ edges and without a copy of $G$ is called an extremal
	graph for $G$. We use $EX(n,G)$ to denote the set of extremal graphs for $G$, i.e.
	\[
	EX(n,G):=\{H\ |\ |V(H)|=n, e(H)= \ex(n,G), \mbox{ and $H$ is $G$-free}\}.
	\]

	The value of $ar(n,G)$ is closely related to the Tur\'an number $\ex(n,G)$ as the following inequality
	\begin{align}\label{low-upp}
		2+\ex(n,\G) \leq ar(n,G) \leq 1+\ex(n,G),
	\end{align}
	where $\G=\{G-e\ |\ e\in E(G)\}$.  Erd\H{o}s, Simonovits and S\'os \cite{Erdos1973}  proved there exists a number $n_0(p)$ such that $ar(n,K_p)=t_{p-1}(n)+2$ for $n> n_0(p)$. Montellano-Ballesteros \cite{Mon-Ball2005} and Neumann-Lara \cite{Mon-Ball2002} extended this result to all values of $n$ and $p$ with $n> p\geq 3$.
	
	An interesting problem concerning anti-Ramsey numbers is to determine the anti-Ramsey number of cycles. Erd\H{o}s, Simonovits and S\'os \cite{Erdos1973} conjectured that $ar(n,C_l) = \left( \frac{l-2}{2} + \frac{1}{l-1} \right) n + O\left( 1 \right)$ and proved it for $l=3$. Alon \cite{Alon1983} proved this conjecture for $l=4$ by showing that $ar(n,C_4)=\lfloor \frac{4n}{3} \rfloor-1$. Jiang, Schiermeyer and West \cite{Jiang} proved the conjecture for $l \leq 7$. Finally, Montellano-Ballesteros and Neumann-Lara \cite{Mon-Ball2005} completely proved this conjecture. Let $kC_3$ denote the union of $k$ vertex-disjoint triangles. Yuan and Zhang \cite{Yuan} provide the exact results of $ar(n,kC_3)$ when $n$ is sufficiently large. Later, Wu and Zhang et al. \cite{Wu} improve the result from $n$ is sufficiently large to $n \geq 2k^2-k+2$. 
	There is a large volume of
	literature on the anti-Ramsey number of graphs.
For those who are interested, I recommend referring to the survey 	by Fujita, Magnant, and Ozeki \cite{Fujita}.
	
	A graph on $2k+1$ vertices consisting of $k$ triangles which intersect in exactly one common vertex is called \textit{friendship graph} (or $ k$-\textit{fan}), denoted by $F_k$.  Erd\H{o}s, F{\"u}redi, Gould and Gunderson \cite{ERDOS1995} determined $ex(n,F_k)$ when $n\geq 50k^2$ and $k\geq 2$. In this paper, our aim is to determine the anti-Ramsey number of friendship graph. 
	
%

\begin{theorem}\label{main1}
	Let $n,k$ be two integers such that $ k\geq 2 $ and $ n\geq 50(k+1)^2 $. The following holds
	\[
	ar(n,F_{k+1})=\ex(n,F_k)+2.
	\]
\end{theorem}

Abbott, Hanson and Sauer~\cite{Abbott} established the value of  $ex(n,\{K_{1,k},kK_2\})$. To finalize the proof of Theorem \ref{main1}, we derive the following result, which determines $ar(n,\{K_{1,k},kK_2\})$ for $k\geq 2, n \geq 3k^2 $.

\begin{theorem}\label{b}
	Let $n,k$ be two integers such that $k\geq 2$ and  $n \geq 3k^2 $. Then
	\begin{equation}
		ar(n,\{K_{1,k+1},(k+1)K_2\})=\begin{cases}
			4, &\mbox{if $ k=2 $},\\
			k^2-k+2, &\mbox{if $ k $ is odd},\\
			k^2-\frac{3}{2}k+2, &\mbox{if $ k $ is even and $ k\ge4 $}.
		\end{cases}	
	\end{equation}
\end{theorem}	

The rest of the paper is organized as follows.
 Section 2 introduces technical lemmas essential for characterizing extremal graphs. Section 3 characterizes $EX(n, \{K_{1,k}, kK_2\})$ and provides the proof of Theorem \ref{b}. Finally, Section 4 completes the proof of Theorem \ref{main1} using the lemmas from Sections 2 and 3.


\section{Some Technical Lemmas on $F_k$ and Matchings}

For a graph $G$, we use $o(G)$ to denote the number of connected components of odd order in $G$. 
\begin{theorem}[Tutte, \cite{Tutte}]\label{Tu47}
A graph $G$ contains a perfect matching if and only if $o(G-T)\leq |T|$ for any $T\subseteq V(G)$.
\end{theorem}

Let $G$ be a graph.  $D_G$ is defined as the set of all vertices in $G$ which are not covered by at least one maximum matching of $G$. The set $A_G$  consists of vertices in $V(G) - D_G$ that are adjacent to at least one vertex in $D_G$.  Finally let $C_G = V(G) - A_G - D_G$. The partition $(D_G,A_G,C_G)$ is called a canonical decomposition. When there is no confusion, we also use $A_G,C_G,D_G$ to denote the induced graphs $G[A_G]$, $G[C_G]$ and $G[D_G]$, respectively.
\begin{theorem}[The Gallai-Edmonds Structure Theorem, \cite{Lovasz1986}]\label{Galli-Edmonds}
For a graph $G$, let $D_G$, $A_G$ and $C_G$ be defined as above. Then the following three statements hold.
\begin{itemize}[itemsep=0pt,parsep=0pt]
	\item [$(i)$] the components of the subgraph induced by $D_G$ are factor-critical;
	\item [$(ii)$] the subgraph induced by $C_G$ has a perfect matching;
	\item [$(iii)$] any maximum matching $M$ of $G$ contains a near perfect matching of each component of $D_G$, a perfect matching of $C_G$ and matches all vertices of $A_G$ with vertices in distinct components of $D_G$.
\end{itemize}
\end{theorem}

Define $ f\left(\nu, \Delta\right)=\max\{e(G)\ |\ \nu(G)\leq \nu, \Delta(G)\leq \Delta\} $. In 1972, Abbott, Hanson and Sauer~\cite{Abbott} determined $ f(k, k)$.
\begin{theorem}[Abbott, Hanson and Sauer,~\cite{Abbott}]\label{kk-value}
Let $k\geq 2$ be an integer. Then
\begin{equation*}
	f(k,k)=\begin{cases}
		k^2+(k-1)/2 ,& \mbox{if k is odd}, \medskip\\
		k^2+k,& \mbox{if k is even}.
	\end{cases}	
\end{equation*}
\end{theorem}
Later Chv\'atal and Hanson\cite{chva1976}  proved the following result, which generalized Theorem \ref{kk-value}.
\begin{theorem}[Chv\'atal and Hanson, \cite{chva1976}]\label{Chvatal}
For every $\nu \geq 1$ and $\Delta\geq 1$,
$$
f( \nu ,\Delta ) =\nu \Delta +\lfloor \frac{\Delta}{2} \rfloor \lfloor \frac{\nu}{\lceil \Delta/2 \rceil} \rfloor \leq \nu \Delta +\nu.
$$
\end{theorem}

Balachandran and Khare \cite{Bala2009} gave a more ``structural" proof of this result. They gave a simple characterization of all the cases where the extremal graph is unique, that is $\nu =1$ or $\lceil \Delta/2 \rceil$ divides $\nu$. For two integers $\nu, \Delta \geq 1$, let $\mathcal{F}_{\nu, \Delta}:=\{G\ |\ \nu(G)\leq \nu, \Delta(G)\leq \Delta, e(G)=f\left(\nu, \Delta\right)\}$. In Section 3, we will give all the extremal graphs for $\nu =\Delta= k-1$ and $\nu = k-2,\Delta=k-1$.

\begin{theorem}[Mantel, \cite{Mantel}]\label{Mantel}
For integer $n\geq 3$, $\ex(n,K_3)=\lfloor n^2/4\rfloor$.
\end{theorem}

Erd\H{o}s et al. \cite{ERDOS1995} proved that $\ex(n,F_2)=\lfloor \frac{n^2}{4} \rfloor +1$ for $n\geq 5$. For $k\geq 3$, they gave the following result.

\begin{theorem}[Erd\H{o}s, F{\"u}redi, Gould and Gunderson, \cite{ERDOS1995}]\label{Fkextremal}
For two integers $n,k$ such that $k\geq 2$ and $n\geq 50k^2 $,   $$
\ex(n,F_k)=\left\{ \begin{matrix}
	\lfloor \frac{n^2}{4} \rfloor +k^2-k, &		\mbox{if $k$ is odd}, \medskip\\
	\lfloor \frac{n^2}{4} \rfloor +k^2-\frac{3}{2}k, &\mbox{if $k$ is even}.\\
\end{matrix} \right.
$$
\end{theorem}

Erd\H{o}s, F{\"u}redi, Gould and Gunderson \cite{ERDOS1995} gave some graph constructions in $EX(n,F_k)$.
Each member in $EX(n,F_k)$ may be obtained from
 a Tur\'{a}n graph $ T_{n,2} $ with a graph $ G $ embedding in one partite set, where $ G $ is empty for $ k=1 $ and $ G\in \mathcal{F}_{k-1,k-1} $ for $ k\geq 2 $.



	\begin{lemma}[Erd\H{o}s, F{\"u}redi, Gould and Gunderson, \cite{ERDOS1995}]	\label{aa}
		If $G$ is $F_{k+1}$-free, then $e(G) \leq \lfloor \frac{n^2}{4} \rfloor +k \Delta(G)$.
	\end{lemma}
	
	
	Let $G$ be a graph with a partiton of the vertices into two non-empty parts $V(G)=V_0\cup V_1$. Define $G_0=G[V_0], G_1=G[V_1]$, and
	$$
	G_{2}=(V(G),\{ v_0 v_1 \in E(G)\ |\  v_0 \in V_0, v_1 \in V_1 \}).
	$$
	
	\begin{lemma}[Erd\H{o}s, F{\"u}redi, Gould and Gunderson, \cite{ERDOS1995}]\label{Gpartition}
		Suppose that $G$ is an $F_{k+1}$-free graph on $ n  $ vertices with $n\geq 24(k+1)^2 $, and with minimum degree $ \delta>(n/2)-(k+1) $, maximum degree $\Delta>n/2 $, Then there exists a partition $V(G)=V_0\cup V_1$, so that $V_0\ne \emptyset$, $V_1\ne \emptyset$, and for each $ i=0,1 $, every $x\in V_i$, the following hold:
		\begin{equation}
			\nu(G_i)\leq k \ \ \ and\ \ \ \  \Delta(G_i)\leq k;
			\label{5}
		\end{equation}
		\begin{equation}
			d_{G_i}(x)+\nu(G_{1-i}[N_G(x)\cap V_{1-i}])\leq k.
			\label{6}
		\end{equation}
	\end{lemma}
	
	\begin{lemma}[Erd\H{o}s, F{\"u}redi, Gould and Gunderson, \cite{ERDOS1995}]\label{y}
		Suppose G is partitioned as above so that (\ref{5})  and  (\ref{6}) are satisfied. If $ G $ is $F_{k+1}$-free, then
		\begin{equation*}
			e(G_0)+e(G_1)-(|V_0|\cdot|V_1|-e(G_{2}))\leq f(k,k).
		\end{equation*}
	\end{lemma}

Let $\mathcal{G}(n,k)$ be the family of graphs where each member is a Tur\'{a}n graph $ T_{n,2} $ with a graph $ G $ embedding in one partite set, where $ G $ is empty for $ k=1 $ and $G\in \mathcal{F}_{k-1,k-1} $ for $ k\geq 2 $.	

\begin{lemma}[Erd\H{o}s, F{\"u}redi, Gould and Gunderson,\cite{ERDOS1995}]\label{extremalgraph}
For $ n\geq 50k^2 $, $EX(n,F_k)=\mathcal{G}(n,k)$.
\end{lemma}		
	
	\section{Anti-Ramsey  Number of $\{K_{1,k},kK_2\}$}
	
	In this section, we determine  $ar(n,\{K_{1,k},kK_2\}$). Basing on Gallai-Edmonds Structure Theorem \cite{Lovasz1986}, we give the detailed description of extremal graphs when $\nu =\Delta= k-1$ and $\nu = k-2$,  $\Delta=k-1$. 
	
Recall that for two integers $\nu, \Delta \geq 1$,  $\mathcal{F}_{\nu, \Delta}:=\{G\ |\ \nu(G)\leq \nu, \Delta(G)\leq \Delta, e(G)=f\left(\nu, \Delta\right)\}$.
	Let $\mathcal{H}_{\nu, \Delta}$ denote the set of graphs $C \cup G(X,Y)$ satisfying the following three conditions:
	\begin{itemize}
		\item  for even $\Delta$, $C$ is a $\Delta$-regular factor-critical graph, and for odd $\Delta$, $C$ is a nearly $\Delta$-regular factor-critical graph;
		\item  $|X|=\nu(G(X,Y))$ and $ \Delta =\Delta(G(X,Y))$ with $ |X|\geq 0 $;
		\item $\nu(G(X,Y))+\frac{|C|-1}{2}=\nu $ and $d_G (v)=\Delta$ for every vertex $v \in X$.
	\end{itemize}
	For odd $\Delta$ and $H:=C\cup G(X,Y)\in \mathcal{H}_{\nu, \Delta}$, let $v'\in V(C)$ such that $d_H(v')=\Delta-1$ and define
	\[
	\Gamma_{\nu, \Delta}(H)=\{(H-uv)+ \{uv'\}\ |\ u\in X, uv\in E(G(X,Y))\}.
	\]
	Let $$
	\mathcal{E}_{\nu, \Delta}=\left\{ \begin{matrix}
		\mathcal{H}_{\nu, \Delta},&		\mbox{if $\Delta$ is even},\medskip\\
		\mathcal{H}_{\nu, \Delta}\cup\bigcup_{H\in \mathcal{H}_{\nu, \Delta}}\Gamma_{\nu, \Delta}(H),&	   	\mbox{if $\Delta$ is odd}.\\
	\end{matrix} \right.
	$$

	%
	%
	
	For odd $k$,
 Balachandran and  Khare~\cite{Bala2009}	characterized $\mathcal{F}(\nu,\Delta)$.
	\begin{lemma}[Balachandran and  Khare,~\cite{Bala2009}]\label{lem99}
		Let $k \geq 3$ be an odd integer. Then $\mathcal{F}_{k-1,k-1}=\{K_k \cup K_k\}$.	
	\end{lemma}
	
	We characterize the case for even $k$.
	\begin{lemma}\label{th99}
		Let $k \geq 4$ be an even integer. Then $\mathcal{F}_{k-1,k-1}=\mathcal{E}_{k-1,k-1}\cup \mathcal{D}_k$, where $\mathcal{D}_k$ denote the set of graphs $K_{k-1}\cup C$
		such that $C$ is a nearly $(k-1)$-regular factor-critical graph of order $k+1$.	
	\end{lemma}
	
	\pf	By the definition of the set $\mathcal{E}_{k-1,k-1} \cup \mathcal{D}_k$, for any graph $G\in \mathcal{E}_{k-1,k-1} \cup \mathcal{D}_k$, it is evident that $\Delta(G)=k-1$  and $\nu(G)=k-1$. Furthermore, the number of edges $e(G)$ is given by the function $f(k-1, k-1)$. Consequently, we can deduce that the set $\mathcal{E}_{k-1,k-1} \cup \mathcal{D}_k\subseteq \mathcal{F}_{k-1,k-1}$.

	Next we show that $\mathcal{F}_{k-1,k-1}\subseteq \mathcal{E}_{k-1,k-1} \cup \mathcal{D}_k$. Let   $H \in \mathcal{F}_{k-1,k-1}$ and let $D_H$, $A_H$ and $C_H$ be defined as previously mentioned. 
	Let $C_1, C_2, \dots, C_t$ denote the connected components of $D_H$. We define $a_H=|A_H|$, $c_H=|C_H|$, and $c_i=|C_i|$ for $i=1,\dots,t$. Without loss of generality, we assume  that  $ c_1\geq c_2\geq \ldots \geq c_t\geq 1$. According to Theorem  \ref{Galli-Edmonds}, $C_H$ has a perfect matching, $C_i$ is factor-critical for $i \in [t]$ and the matching number of $H$ is
	\begin{align}\label{GE-equality}
		k-1=a_H+\frac{c_H}{2}+\sum\limits_{i=1}^t{\frac{c_i-1}{2}}.
	\end{align}
	
	\medskip
	\textbf{Claim 1.}~For each $v \in A_H$, $d_{A_H}(v)=0$ and $d_H(v)=k-1$.
	\medskip
	
	If there is an edge $e=uv$ with $u,v \in A_H$, we may obtain a graph $H'$ with $e(H)+1$ edges from $H$ by deleting edge $e$ and adding a new vertex with joining it to both $u$ and $v$. By Theorem \ref{Galli-Edmonds} (iii), we have $\nu(H') = k-1$, $\Delta(H') = k-1$ and $e(H')=e(H)+1$, a contradiction since $H\in  \mathcal{F}_{k-1,k-1}$. If there is a vertex $u \in A_H$ such that $d_H(u)<k-1$. We obtain a graph $H''$ from $H$ by adding a new vertex and joining it to $u$. By Theorem \ref{Galli-Edmonds} (iii), $\nu(H'') = k-1$, $\Delta(H'') = k-1$ and $e(H'')=e(H)+1$, which contradicts the fact that $H\in \mathcal{F}_{k-1,k-1}$. This completes the proof of Claim 1.
	
	\medskip
	\textbf{Claim 2.}~$C_i$ is complete or nearly $(k-1)$-regular for $i \in [t]$.
	\medskip

Firstly, let's assume that  $c_i \leq k-1$. We assert that $C_i$ is a complete graph. If not, suppose there are two non-adjacent vertices $u,v\in V(C_i)$. We add two new vertices $u',v'$  and define a new graph $H'$ with the vertex set $V(H)\cup \{v',u'\}$ and edge set
	\[
	(E(H)-E_H(\{u,v\},A_H))\cup \{uv\}\cup \{u'x\ |\ x\in N_H(u)\cap A_H\}\cup \{v'x\ |\ x\in N_H(v)\cap A_H\}.
	\]
	One can observe that  $\nu(H')\leq k-1$, $\Delta(H')=k-1$ and $e(H')=e(H)+1$ contradicting the choice of $H$. Moving forward, we consider the case where $c_i\geq k+1$.
We claim that $C_i$  is nearly $(k-1)$-regular.
If not, suppose $C_i$ is not nearly \((k-1)\)-regular. Given that $k\geq 4$, there exists a nearly $(k-1)$-regular factor-critical graph of order $|C_i|$, denoted by $C_i'$.
	Let $T:=\{y\in V(C_i)\ |\ y\notin N_H(A_H)\}$, and let $H'$ be a graph with vertex set
	$
	(V(H)-T)\cup V(C_i')
	$
	and edge set
	\[
	E(H-V(C_i))\cup E_H(A_H,V(C_i)-T)\cup E(C_i').
	\]
	Since $e(C_i')>e(C_i)$, we have $e(H')>e(H)$.  According to Theorem \ref{Galli-Edmonds} (iii), we find that $\nu(H') = k-1$ and $\Delta(H') = k-1$. Consequently, we deduce that
$H'\in \mathcal{F}_{k-1,k-1}$, which contradicts the fact that $H\in \mathcal{F}_{k-1,k-1}$. This completes the proof of Claim 2.
	
	\medskip
	\textbf{Claim 3.}~$c_2 =k-1$ or $c_2 =1$.
	
	\medskip
	
	By contradiction, suppose that $c_2 \notin \{1,k-1\}$. According to equation (\ref{GE-equality}), $\frac{c_1-1}{2}+ \frac{c_2-1}{2} \leq k-1$. Thus we have $c_1+c_2 \leq 2k$, which implies $c_2\leq k$.
	Since $c_2$ is odd and $k$ is even, it follows that $3\leq c_2 \leq k-3$. Define $T$ as the set of vertices in $C_1\cup C_2$ that are not in the neighborhood of $A_H$, i.e.,
 $T:=V(C_1\cup C_2)-N_H(A_H)$ and let $T':=V(C_1\cup C_2)-T$.
	
	Firstly,  consider the case where $c_1+c_2\geq k+2$. Then there exists a nearly $(k-1)$-regular graph with $c_1+c_2-1$ vertices denoted by $C_1'$.
	Now we can construct a new graph $H'$ with vertex set $(V(H)-T)\cup V(C_1')$ and edge set
	\[
	E(H-V(C_1\cup C_2))\cup E_H(A_H,T')\cup E(C_1').
	\]
	It is clear that that $\nu(H') \leq k-1$, $\Delta(H') = k-1$ and
	\begin{align*}
		e(H')-e(H)&\geq  \frac{(c_1+c_2-1)(k-1)-1}{2} - \frac{c_1(k-1)-1}{2} - \frac{c_2(c_2-1)}{2} \\
		&= \frac{(c_2-1)(k-1-c_2)}{2}>0,
	\end{align*}
	which contradicts the fact that $H\in \mathcal{F}_{k-1,k-1}$.

	Secondly,  consider the case where $c_2 + c_1 \leq k$. Recall that $c_2\geq 3$. Then we have $c_1 \leq k-3$ and $k\geq 6$. Now  we may construct a  graph $H'$
	with vertex set $(V(H)-T)\cup V(K_{c_1+c_2-1})$ and edge set
	\[
	E(H-V(C_1\cup C_2))\cup E_H(A_H,T')\cup E(K_{c_1+c_2-1}).
	\]
	Similarly, one can see that $\nu(H') \leq k-1$, $\Delta(H') = k-1$ and
	\begin{align*}
		e(H')-e(H)&\geq {c_1+c_2-1\choose 2}-{c_1\choose 2}-{c_2\choose 2}\\
		&=(c_1-1)(c_2-1)>0,
	\end{align*}
which contradicts the fact that $H\in \mathcal{F}_{k-1,k-1}$. This completes the proof of Claim 3.

\medskip
\textbf{Claim 4.}~$c_H=0$.

\medskip

By contradiction, suppose that $c_H \geq 2$. By (\ref{GE-equality}), we have $c_H+c_1 \leq 2k-1$. So either $c_H \leq k-2$ or $c_1 \leq k-1$ holds. Let $T:=(C_H\cup V(C_1))-N_H(A_H)$ and $T':=(C_H\cup V(C_1))-T$.


Consider $c_H +c_1 \geq k+1$. Let $C_1'$ be  a nearly $(k-1)$-regular factor-critical graph with $c_H +c_1$ vertices. Then we may construct a graph $H'$ with vertex set $(V(H)-T)\cup V(C_1')$ and edge set
\[
E(H-(V(C_1)\cup C_H))\cup E(C_1')\cup E_H(A_H,T').
\]
Note that
\begin{align*}
e(C_1')&=\frac{1}{2}((c_H+c_1)(k-1)-1)\\
&= \frac{1}{2}c_H(k-1)+\frac{1}{2}(c_1(k-1)-1)\\
&>e(C_H)+e(C_1) \quad\mbox{(since $c_H+c_1\leq 2k-1$)}.
\end{align*}
Indeed, one can observe  that $\nu(H') \leq k-1$, $\Delta(H') = k-1$ and
$e(H') > e(H)$, which contradicts the definition of $H$.

Now we may assume that $c_H +c_1 \leq k-1$. We construct a graph $H'$ with vertex set $(V(H)-T)\cup V(K_{c_1+c_H})$ and edge set
\[
E(H-V(C_1\cup C_H))\cup E_H(A_H,T')\cup E(K_{c_1+c_H}).
\]
One can see that $\nu(H') \leq k-1$, $\Delta(H') \leq k-1$ and
\[
e(H')-e(H) \geq {c_1+c_H\choose 2}-{c_1\choose 2}-{c_H\choose 2}=c_1c_H>0,
\]
which contradicts the fact that $H\in \mathcal{F}_{k-1,k-1}$. This completes the proof of Claim 4.
\medskip

By Claim 3, we have $c_2\in \{1,k-1\}$. Firstly, consider the case where $c_2=k-1$. By Claim 2, $C_2$ is a complete graph. Recall that $c_1+c_2\leq 2k.$ It follows that $c_1\in \{k-1,k+1\}$.  According to Claim 4, we have $c_H=0$. If $c_1=k-1$,  equation (\ref{GE-equality}) yields two possibilities: $a_H=0$ and $c_3 = 3$ or $a_H=1$ and $c_3 = 1$. So we have
\begin{align*}
	e(H)-f(k-1,k-1)&=e(H)-{k-1\choose 2}-(k^2-2)/2\\
	&\leq {k-1\choose 2}+(k-1)-(k^2-2)/2=-k/2 + 1<0,
\end{align*}
a contradiction. Thus we may assume that $c_1=k+1$, which implies that $a_H+c_H=0$. It follows that $ e(C_1)=(k^2-2)/2$ since $c_2=k-1$. Note that $C_1$ is factor-critical. Since $\Delta(C_1)\leq k-1$, $C_1$ is a nearly $(k-1)$-regular factor-critical graph.
Consequently,  $H\in \mathcal{D}_k$.

Next we consider the case where  $c_2=1$. We claim $c_1\geq k+1$. Otherwise, suppose $c_1\leq k-1$. Then by (\ref{GE-equality}) and Claims 3 and 4, we have
\[
e(H)= {c_1\choose 2}+a_H(k-1)\leq (k-1)^2<f(k-1,k-1),
\]
a contradiction. By Claim 2, $C_1$ is a nearly $(k-1)$-regular factor-critical graph, which implies that $e_H(V(C_1),A_H)\leq 1$ since $\Delta(H)\leq k-1$. Thus we have $H\in \mathcal{E}_{k-1,k-1}$.
This completes the proof of Lemma \ref{th99}.   \qed
%
%

\begin{lemma}\label{th10}
	For $k \geq 3$, $\mathcal{F}_{k-2, k-1}=\mathcal{E}_{k-2, k-1}$.
\end{lemma}

\pf According to the definition of the set $\mathcal{E}_{k-2, k-1}$, for any $H\in \mathcal{E}_{k-2, k-1}$, it is clear that $\Delta(H)=k-2$, $\nu(H)=k-1$ and $e(H)=f(k-2,k-1)$. Therefore, we can conclude that $\mathcal{E}_{k-2, k-1}\subseteq \mathcal{F}_{k-2, k-1}$.

Next we show that $\mathcal{F}_{k-2, k-1}\subseteq \mathcal{E}_{k-2, k-1}$. Suppose that $H \in \mathcal{F}_{k-2, k-1}$. Let $D_H$, $A_H$ and $C_H$ be defined  as previously mentioned. Let $C_1, C_2, \ldots, C_t$ be the connected components of $D_H$. Define  $c_H=|C_H|$, $a_H=|A_H|$, $c_i=|C_i|$ for $i=1,\dots,t$. Without loss of generality, we assume that $ c_1\geq c_2\geq \ldots\geq c_t\geq 1 $. By Theorem \ref{Galli-Edmonds}, we have
\begin{align}\label{GE-eq2}
	k-2=a_H+\frac{c_H}{2}+\sum\limits_{i=1}^t{\frac{c_i-1}{2}}.
\end{align}

Following the same discussion as in Claims 1 and 2, we arrive at the following claim.

\medskip
\textbf{Claim 5.}~For each $v \in A_H$, $d_{A_H}(v)=0$ and $d_H(v)=k-1$.  Each connected component $C_i$ of $D_H$ is factor-critical graph and  it satisfies one of the following conditions:
\begin{itemize}
	\item [$(a)$] $C_i$ is a complete graph;
	\item [$(b)$] if $k$ ie even, $C_i$ is  a nearly $(k-1)$-regular graph;
	\item [$(c)$] if $k$ is odd, $C_i$  is  a $(k-1)$-regular graph.
\end{itemize}
\medskip

%
%

\medskip
\textbf{Claim 6.}~$c_2 =1$.
\medskip

By contradiction, suppose that $c_2 \geq 3$. From (\ref{GE-eq2}),  it follows that $\frac{c_1-1}{2}+ \frac{c_2-1}{2} \leq k-2$, hence we have $c_1+c_2 \leq 2k-2$.
Since  $c_1+c_2$ is even, we assert that  $c_1+c_2\leq 2\lfloor k/2\rfloor$.
If not, then  $c_1+c_2\geq 2\lfloor k/2\rfloor+1$, leading to $c_1+c_2\geq k+2$ for even $k$ and $c_1+c_2\geq k+1$ for odd $k$. Let $T:=V(C_1\cup C_2)-N_H(A_H)$ and $T':=V(C_1\cup C_2)-T$. For even $k$, let $C_1'$ be a  nearly $(k-1)$-regular factor-critical graph
with $c_1+c_2-1$ vertices; for odd $k$, let $C_1'$ be a $(k-1)$-regular factor-critical graph with $c_1+c_2-1$ vertices.
We construct a graph $H'$ with vertex set $(V(H)-T)\cup V(C_1')$ and edge set
\[
E(H-V(C_1\cup C_2))\cup E(C_1')\cup E_H(A_H,T').
\]
Note that $\nu(H') \leq k-2$, $\Delta(H') = k-1$. If $c_2<k-1$, then
\begin{eqnarray*}
	e(H')-e(H)&\geq& \lfloor \frac{(c_1+c_2-1)(k-1)}{2} \rfloor-\lfloor \frac{c_1(k-1)}{2} \rfloor- \frac{c_2(c_2-1)}{2} \\
	&=& \frac{(c_2-1)(k-1-c_2)}{2}>0,
\end{eqnarray*}
which contradicts the fact  that $H \in \mathcal{F}_{k-2, k-1}$. So we may assume that $c_1=c_2=k-1$ and $k$ is even. Then
\begin{align*}
	e(H')-e(H)& = \frac{(2k-3)(k-1)-1}{2}-2{k-1\choose 2} \\
	&= \frac{k-1}{2}>0,
\end{align*}
a contradiction since $H \in \mathcal{F}_{k-2, k-1}$. So we may assume that $c_1+c_2\leq k$ for even $k$ and $c_1+c_2\leq k-1$ for odd $k$. Let $H'$ be a graph with vertex set $(V(H)-T)\cup V(K_{c_1+c_2-1})$ and edge set
\[
E(H-V(C_1\cup C_2))\cup E_H(A_H,T')\cup E(K_{c_1+c_2-1}).
\]
Note that $\nu(H')\leq k-2$ and $\Delta(H')\leq k-1$
still holds. Moreover, we have
\begin{eqnarray*}
	e(H')-e(H)&=& \left( \begin{array}{c}
		c_1+c_2-1\\
		2\\
	\end{array} \right) -\left( \begin{array}{c}
		c_1\\
		2\\
	\end{array} \right) -\left( \begin{array}{c}
		c_2\\
		2\\
	\end{array} \right) \\
	&=&(c_1-1)(c_2-1) >0 ,
\end{eqnarray*}
a contradiction since $H \in \mathcal{F}_{k-2, k-1}$. This completes the proof of Claim 6.

\medskip
\textbf{Claim 7.}~$c_H=0$.

\medskip

By contradiction, suppose that $c_H \geq 2$. By equation (\ref{GE-eq2}),  $c_H+c_1 \leq 2k-3$. Therefore either $c_H\leq k-2$ or $c_1\leq k-2$ holds.   Let $T:=V(C_1\cup C_H)-N_H(A_H)$ and $T':=V(C_1\cup C_H)-T$.
We claim $c_1+c_H\leq k-1$.  	
Otherwise, suppose $c_H +c_1 \geq k$. For odd $k$, let $C_1'$ be a $(k-1)$-regular factor-critical graph of order $c_1+c_H$
and for even $k$, let $C_1'$ be a nearly $(k-1)$-regular factor-critical graph of order $c_1+c_H$.
Now, we construct a graph $H'$ from $H$ with vertex set $(V(H)-T)\cup V(C_1')$ and edge set
\[
E(H-V(C_1\cup C_H))\cup E_H(A_H,T')\cup E(C_1').
\]
One can see that $\nu(H') \leq k-2$, $\Delta(H') = k-1$. Recall that $\min \{c_H,c_1\}\leq k-2$. Note that
\begin{align*}
	e(H')-e(H)&=\lfloor(c_1+c_H)(k-1)/2\rfloor-e(C_H)-e(C_1)\\
	&\geq \lfloor(c_1+c_H)(k-1)/2\rfloor -\max\{c_H,c_1\}(k-1)/2-{\min\{c_1,c_H\}\choose 2}\\
	&>0,
\end{align*}
a contradiction since $H\in \mathcal{F}_{k-2, k-1}$.
Therefore, we may assume that  $c_1+c_H\leq k-1$. We construct a graph $H'$ with  vertex set $(V(H)-T)\cup V(K_{c_1+c_H})$ and edge set
\[
E(H-V(C_1\cup C_H))\cup E_H(A_H,T')\cup E(K_{c_1+c_H}).
\]
Clearly, $\nu(H') \leq k-2$, $\Delta(H') = k-1$ and
\begin{align*}
	e(H')-e(H) \geq {c_1+c_H\choose 2}-{c_1\choose 2}-{c_H\choose 2}=c_1c_H>0,
\end{align*}
a contradiction again. This completes the proof of Claim 7.
\medskip

By  Claims 5, 6 and 7 and equation (\ref{GE-eq2}), we have $k-2=a_H+\frac{c_1-1}{2}$. Thus,  we can derive
\begin{equation*}
	f(k-2,k-1)=\begin{cases}
		a_H(k-1)+\frac{c_1(k-1)}{2}, & \mbox{if $ k $ is odd}, \medskip\\
		a_H(k-1)+\frac{c_1(k-1)-1}{2}, & \mbox{if $ k $ is even}.
	\end{cases}	
\end{equation*}
Hence, when $k$ is odd, $C_1$ is a $(k-1)$-regular graph, and there is no edge between $C_1$ and $A_H$. When $k$ is even, $C_1$ is a nearly $(k-1)$-regular graph with a vertex $u$ satisfying $d_{C_1}(u)=k-2$ and $E_H(\{ u \},A_H)\leq 1$. 	
This completes the proof.    \qed

\begin{lemma}\label{kk-anti-lower}
	For $k\geq 2, n \geq 3k $,
	\begin{equation}
		ar(n,\{K_{1,k+1},(k+1)K_2\})\geq\begin{cases}
			4, & \mbox{if  $k=2$},\\
			k^2-k+2, & \mbox{if $k$ is odd},\\
			k^2-\frac{3}{2}k+2, & \mbox{if $k$ is even and $k \geq 4$}.
		\end{cases}	
	\end{equation}
	
\end{lemma}	
\pf Let's consider the case when $k=2$.  We denote the vertex set of the complete graph $K_n$ as $V(K_n)=\{x_1,\ldots,x_n\}$.  Define a function  $f:E(K_n)\rightarrow [3]$ such that
\begin{equation}
	f(x_ix_j)=\begin{cases}
		1, & \mbox{if $i=1$ and $j=2$},\\
		2, & \mbox{if $i=1$ and $3\leq j\leq n$},\\
		3,& \mbox{otherwise}.
	\end{cases}	
\end{equation}
It is evident that the colored graph does not contain a rainbow $K_{1,3}$ or $3K_2$. Therefore, we have  $ar(n,\left\{K_{1,3},3K_2\right\})\geq 4$.

For odd $k\geq 3$, let $K_k^1$ and $K_k^2$ be two vertex-disjoint subgraphs of $K_n$. Define a bijective coloring $h:E(K_k^1\cup K_k^2)\rightarrow [k^2-k]$.
Let  $G_1$ be the $n$-vertex complete graph $K_n$ with edge-coloring $h_{G_1}$, where
\begin{equation}
	h_{G_1}(e)=\begin{cases}
		h(e), & \mbox{if $e\in E(K_k^1\cup K_k^2)$},\\
		0,& \mbox{otherwise}.
	\end{cases}	
\end{equation}
Since the matching number of $K_{k}^1\cup K_k^2$ is $k-1$, then $G_1$ contains no rainbow matchings of size $k+1$.
It is clear that $G_1$ contains no rainbow   $K_{1,k+1}$. Thus for odd $k$, $ar(n,\{K_{1,k+1},(k+1)K_2\})\geq k^2-k+2$.

For even $k\geq 4$, let $C$ be a nearly $(k-1)$-regular factor-critical graph with order $k+1$.
We select two vertex-disjoint subgraphs of $K_n$, saying  $C'$ and $K_{k-1}$  such that $C'$ is a copy of $C$.
Define a bijective coloring  $g:E(K_{k-1}\cup C')\rightarrow [k^2-\frac{3}{2}k]$. Let  $G_2$ be the $n$-vertex complete graph $K_n$ with edge-coloring $g_{G_2}$, where
\begin{equation}
	g_{G_2}(e)=\begin{cases}
		g(e), & \mbox{if $e\in E(K_{k-1}\cup C')$},\\
		0,& \mbox{otherwise}.
	\end{cases}	
\end{equation}
Since $\nu(K_{k-1}\cup C')=k-1$, then $G_2$ contains no rainbow matchings of size $k+1$. Clearly, $G_2$ contains no rainbow $K_{1,k+1}$.
Thus for even $k\geq 4$, $ar(n,\{K_{1,k+1},(k+1)K_2\})\geq k^2-3k/2+2$. This completes the proof.
\qed

\begin{lemma}\label{edge-vertex-critical}
	Let $r\geq 4$ be an integer and $G$ be a nearly $r$-regular graph such that $|G|\leq 2r-1$ and $|G|$ is odd. Then for any $v\in E(G)$ and any $e\in E(G)$, $G-v-e$ has a perfect matching.
\end{lemma}

\pf Suppose that the result does not hold. Then there exists $v\in V(G),e\in E(G)$ such that $G-v-e$ contains no perfect matching. Let's denote $G-v-e$ by $F$. We claim $G$ contains no nontrivial edge-cut of size no more than $r-1$. Otherwise, let $Q$ be a nontrivial minimal edge-cut of size no more than $r-1$ in $G$. By the minimality of $Q$, $G-Q$ consists of two connected components, denoted by $G_1$ and $G_2$ such that $2\leq |G_1|< |G_2|$. Then we have
\begin{align*}
	r-1\geq |Q|&= \sum_{x\in V(G_1)}(d_G(x)-d_{G_1}(x))\\
	&\geq (r-(|G_1|-1))|G_1|-1\\
	&\geq 2r-3,
\end{align*}
which leads to a contradiction given that $r\geq 4$. By Theorem \ref{Tu47}, there exists a subset $T\subseteq V(G)-v$ such that $q=o(F-T)\geq |T|+2$. We use $C_1,\ldots, C_q$ to denote these odd components.
Then we  may infer that
\begin{align*}
	r(|T|+1)&= \sum_{x\in T\cup\{v\}}d_G(x)\\
	&\geq \sum_{i=1}^q e_G(V(C_i),T\cup \{v\})-2\\
	&\geq qr-3\geq r(|T|+2)-3,
\end{align*}
a contradiction since $r\geq 4$.
This completes the proof. \qed

\medskip

\noindent\textbf{Proof of Theorem \ref{b}.} Let \begin{equation}
	c(n,k)=\begin{cases}
		4, &\mbox{if $ k=2 $},\\
		k^2-k+2, &\mbox{if $ k $ is odd},\\
		k^2-\frac{3}{2}k+2, &\mbox{if $ k $ is even and $ k\ge4 $}.
	\end{cases}	
\end{equation}
According to Lemma \ref{kk-anti-lower},  we have established that $ar(n,\{K_{1,k+1},(k+1)K_2\})\geq c(n,k)$.  Therefore, it suffices to prove that  $ar(n,\{K_{1,k+1},(k+1)K_2\})\leq c(n,k)$. Let $c:E(K_n)\rightarrow [c(n,k)]$ be a surjective  edge-coloring of $K_n$. The
edge-colored $K_n$ is represented  by $H$. By way of contradiction, assume that $H$ contains no neither rainbow $K_{1,k+1}$ nor rainbow $(k+1)K_2$.
Let $G$ be a rainbow edge-induced subgraph of $H$ with exactly $c(n,k)$ edges. We can assume
that the maximum degree  $\Delta(G)\leq k$ and the matching number $\nu(G)\leq k$. Given that $n\geq 3k^2>2c(n,k)+2$,  there exist two vertices $x,y\in V(H)-V(G)$. We remove the edge with color $c(xy)$ from $G$ and denote the resulting graph as $G'$. Define $V_k:=\{x\in V(G')\ |\ d_{G'}(x)=k\}$. Without loss of generality, we may select $G$ and $G'$ such that $|V_k|$ is minimal.

Let's consider the case when $k=2$. Given that $\Delta(G)\leq 2$, $\nu(G)\leq 2$ and $e(G)=4$, the graph $G$ must be either $ K_3\cup K_2 $, a path of length four or
consists of two vertex-disjoint paths of length two. If $G$ is $ K_3\cup K_2 $, we assume that the edge  set  $ E(K_3\cup K_2)=\{v_1v_2,v_2v_3,v_1v_3, v_4v_5\}$. Since $ H $ contain no rainbow $ K_{1,3} $,  $ c(v_1w)\in \{c(v_1v_2),c(v_1v_3)\} $ for any $ w\in V(G)-\{x,y\}\cup \{v_i\ |\ i\in [5]\} $. This implies that $ P_4 =wv_1v_2v_3$ or $ P_4 =wv_1v_3v_2$ is a rainbow path in $ H $. In any case, we may find a matching $M$ of size three  in $H$, which is a contradiction.

Now, let's discuss the case when $ k\geq 3 $.  If $G'$ has a matching $M_0$ of size $k$, then $M_0\cup \{xy\}$ would be a rainbow matching of size $k+1$ in $H$. So we may assume that $\nu(G')\leq k-1$.

\medskip
\textbf{Claim 8.}~$V_k\neq \emptyset$ and $G'[V_k]$ is a clique of $G'$.
\medskip

Since $ e(G')=c(n,k)-1>f(k-1,k-1)$, we can infer that $\Delta(G') = k$. Therefore, $V_k\neq \emptyset$.		
Next, we show that $G'[V_k]$ is a clique by contradiction. Assume that there exist two vertices $u,v\in V_k $ such that $uv\notin E(G')$. Consider the graph $G'+uv$.
In this case,  either the set $\{uw\ |\ w\in N_{G'}(u)\}\cup \{uv\}$ or the set $\{vw\ |\ w\in N_{G'}(v)\}\cup \{uv\}$ induces a rainbow $K_{1,k+1}$ in $H$, which is  a contradiction.
This completes the proof of Claim 8.

%
Let $V_k :=\{x_1,\ldots, x_r\}$, where $1\leq r\leq k+1$. Define $ S:= V(H) - V(G')$. Note that $|S| \geq n-2c(n,k)\geq k^2$.
For every $x\in V(H)$ and $Q\subseteq V(H)$, let $\mathcal{N}_Q^c(x):=\{c(xu)\ |\ u\in Q\}$. when $Q=V(H)$, we simply write  $\mathcal{N}_Q^c(x)$ for $\mathcal{N}^c(x)$ .

\medskip
\textbf{Claim 9.}~For any $uv\in E(H)$ with $u\in V_k$ and $v\in S$, there exists an edge $uv'\in E(G')$ such that $v'\notin V_k$ and $c(uv')=c(uv)$.
\medskip

If $c(uv)\notin \{c(uw)\ |\ w\in N_{G'}(u)\}$, then  it's clear that adding the edge $uv$ to $G'$ would result in a rainbow
 $K_{1,k+1}$, which is a contradiction. Therefore, we can assume that $c(uv)\in \{c(uw)\ |\ w\in N_{G'}(u)\}$. Suppose there exists an edge $uv'$ such that $c(uv)=c(uv')$ and $v'\in V_k$. Note that $c(uv')\neq c(xy)$. Then, the number of vertices with degree $k$ in the graph
  $(G'-uv')+uv$ is less than $|V_k|$, which contradicts the selection of $G'$. This completes the proof of Claim 9.

\medskip
\textbf{Claim 10.}~There is a rainbow matching $M$ of size $r$  in $E_H(V_k,S-\{x,y\})$ and $r\leq k-2$.
\medskip



%

Recall that $G'$ is a rainbow subgraph of $H$. By Claim 9, for any two distinct vertices $u,v\in V_k$, we have $\mathcal{N}^c_S(u)\cap \mathcal{N}^c_S(v)=\emptyset$. Therefore, $E_H(V_k,S-\{x,y\})$ contains a rainbow matching $M$ of size $r$. This implies that $r=|V_k|\leq k$.

By Claim 8, we have
\begin{align*}
	e(G'-V_k)&= c(n,k)-1-{|V_k|\choose 2}-(k+1-|V_k|)|V_k|\\
	&= c(n,k)-1+|V_k|^2/2-(k+1/2)|V_k|.
\end{align*}
By Claim 9, every edge in $G'-V_k$ has different colors from all edges of $M$. If $|V_k|=k$, then $e(G'-V_k)=c(n,k)-1-k(k+1)/2>0$. Therefore, by Claim 9, for any  $e\in E(G'-V_k)$, $\{e\}\cup M$ is a rainbow matching of size $k+1$, which is a contradiction. Hence, it suffices to show that $r\neq k-1$.

Suppose that $r=k-1$. Then we have
\begin{equation}\label{eq12}
	e(G'-V_k)\geq\begin{cases}
		k(k-3)/2+2, & \mbox{if $k$ is odd},\\
		k(k-4)/2+2, & \mbox{otherwise}.
	\end{cases}	
\end{equation}
 According to equation (\ref{eq12}),  $e(G'-V_k)>0$ for $k\geq 3$. Let $e\in E(G'-V_k)$. By Claim 9, $c(e)\notin \{c(e')\ |\ e'\in M\}$. Recall that $G'+xy$ is a rainbow subgraph of $H$. Therefore, $M\cup \{xy,e\}$ forms a rainbow matching of size $k+1$ in $H$, which is a contradiction. This completes the proof of Claim 10.


%

By the definition of $V_k$ and by Claim 10,  we have $ \Delta(G'-V_k)\leq k-1 $ and  $\nu(G'-V_k)\leq k-1-r$. If $\Delta(G'-V_k)\leq k-2$, then by Theorem \ref{Chvatal},
\begin{align*}
	e(G')&= e(G'-V_k)+{r\choose 2} +r(k-r+1)\\
	&\leq f(k-1-r,k-2)+kr-{r\choose 2}\\
	&\leq (k-1)(k-1-r)+kr-{r\choose 2}\\
	&=k^2-2k+1+r-{r\choose 2}\\
	&\leq k^2-3k/2+1/2\\
	&<f(k-1,k-1)+1,
\end{align*}
a contradiction. So we may assume that $\Delta(G'-V_k)=k-1$. Recall that $r\geq 1$.   One can see that
\begin{eqnarray*}
	e(G')&\leq &f(k-1-r,k-1)+kr-{r\choose 2}\\
    &\leq& (k-1-r)(k-1)+\lfloor(k-1)/2\rfloor \lfloor \frac{k-1-r}{\lceil (k-1)/2\rceil}\rfloor+kr-{r\choose 2}\quad \mbox{(by Theorem \ref{Chvatal})}\\
	&\leq&(k-1-r)(k-1)+\lfloor(k-1)/2\rfloor+kr-{r\choose 2}\quad \mbox{(since $r\geq 1$)}\\
	&=&(k-1)^2+\lfloor(k-1)/2\rfloor-\frac{1}{2}r(r-3)\\
	&\leq &\begin{cases}
		k^2-k, &\mbox{if $ k $ is odd},\medskip\\
		k^2-\dfrac{3}{2}k+r-{r\choose 2},  & \mbox{if $ k $ is even}.
	\end{cases}	
\end{eqnarray*}
For odd $k$, we see that  $e(G') < k^2-k+1$, which is a contradiction.  When $k$ is even, we have $e(G') < k^2-\frac{3}{2}k+1$ for $r\geq 3$, a contradiction. Therefore, we may assume that $r\in \{1,2\}$ and $k$ is even. Moreover, we can infer that
\begin{align*}
	f(k-1,k-1)+1\leq e(G')\leq k^2-\dfrac{3}{2}k+1,
\end{align*}
i.e.,
\begin{align}
	f(k-1,k-1)+1=e(G')= k^2-\dfrac{3}{2}k+1.
\end{align}
Since $r\leq 2$, and $G'[V_k]$ is a clique, there exists an edge $f=uv$ such that $\Delta(G'-f)\leq k-1$. Recall that $\nu(G'-f)\leq \nu(G')\leq k-1$.
  Let $F$ be the  edge-induced subgraph of $G'$ induced by $E(G')-f$. Thus, by Lemma \ref{th99}, $F\in \mathcal{F}_{k-1,k-1}$. Without loss of generality, assume that $d_{G'}(u)=k$. We will now consider three cases.

\medskip
\textbf{Case 1.}~$F$ is a nearly $(k-1)$-regular factor-critical graph of order $2k-1$.
\medskip

If $|\{u,v\}\cap V(F)|\leq 1$, then $G'-\{u,v\}$ contains a matching $M_1$ of size $k-1$. It follows that $M_1\cup \{f,xy\}$ is a rainbow matching of size $k+1$ in $H$, which is a contradiction. Therefore, we may assume  $\{u,v\}\subseteq V(F)$. Since $H$ contains no rainbow $K_{1,k+1}$, for any $uw$ with $w\in S-\{x,y\}$, $c(uw)\in \{c(uw')\ |\ w'\in N_{G'}(u)\}$. Let $M_2$ be a perfect matching of $F-u$. Then for $w\in S-\{x,y\}$, $\{xy, uw\}\cup M_2$ is a rainbow matching of size $k+1$ in $H$, which is a contradiction.

\medskip
\textbf{Case 2.}~$F$ is the union of $K_{k-1}$ and graph $C$, where $C$ is a nearly $(k-1)$-regular factor-critical graph with order $k+1$.
\medskip

	
Recall that $d_{G'}(u)=k$. Therefore, we have $u\in  V(C)$. Since $H$ contains no rainbow $K_{1,k+1}$, then for any $uw$ with  $w\in S-\{x,y\}$, the color $c(uw)$  must be in the set
 $\{c(uw')\ |\ w'\in N_{G'}(u)\}$. Let $M_3$ be a perfect matching of  $C-u$ and let $M_4$ be a near perfect matching of $K_{k-1}$ that avoids the color $c(uw)$. Then, the union $M_3\cup M_4\cup \{xy,uw\}$ forms a rainbow matching of size $k+1$ in $H$, which is a contradiction.

\medskip
\textbf{Case 3.}~There exists a subset $\mathcal{E}_0\subseteq E(F)$  such that $|\mathcal{E}_0|\leq 1$ and  $F-\mathcal{E}_0$ is the union of $C$ and bipartite graph $R$ with bipartition $(X,Y)$.
\medskip

By the definition  of $\mathcal{F}_{k-1,k-1}$, $F$ satisfies  the following conditions.
\begin{itemize}
	\item [$(i)$] $C$ is a nearly $(k-1)$-regular factor-critical graph of order at least $k+1$;
	\item [$(ii)$] $|X|+(|C|-1)/2=k-1$ and $|X| \geq 1$;
	\item [$(iii)$] $d_{F}(w)=k-1$ for all $w\in X$.
\end{itemize}
Since $|C|\geq k+1$ and $|X|\geq 1$,  by Condition (ii),  we have $|C|\leq 2k-3$ and $|X|\leq k/2-1\leq k-3$. For any edge  $h\in E(R)$, by (iii), we may greedily find a matching  of size $|X|$ in $R-V(C)-h$.
 Let $U_{k-1}\subseteq V(C)$ denote the vertex subset of degree $k-1$ in $F$. Moreover,  one can see that $u\in U_{k-1}\cup X$ since $|X|\leq k/2-1$, $d_{F}(u)=k-1$ and $R$ is a bipartite graph.

Consider $u\in U_{k-1}$. Let $w\in S-\{x,y\}$. Since $H$ contains no rainbow $K_{1,k+1}$, we have $c(uw)\in \{c(uw')\ |\ w'\in N_{G'}(u)\}$, which implies that $c(uw)\neq c(xy)$.  Let $M_5$ be a matching of size $|X|$ in $R$ and let $M_6$ be a perfect matching of $C-u$. Then, $M_5\cup M_6\cup \{xy,uw\}$ is a rainbow matching of size $k+1$, a contradiction.

Next, we may assume that $u\in X$.  Since $|\mathcal{E}_0|\leq 1$ and $k\geq 3$, we can choose one vertex $v_1\in U_{k-1}$ such that $uv_1\notin E(G')$. Recall $H$ contains no rainbow $K_{1,k+1}$, which implies that $c(uv_1)\in \{c(uw')\ |\ w'\in N_{G'}(u)\}$. Let $v_2\in S-\{x,y\}$. With a similar discussion, we have $c(v_1v_2)\in \{c(uv_1)\}\cup \{c(v_1w')\ |\ w'\in N_{G'}(v_1)\}$, otherwise, $H$ would  contain a rainbow  $K_{1,k+1}$ with center $v_1$, a contradiction. Let $M_1'$ be  a matching of size $|X|$ in $R$ avoiding the edge colored by $c(v_1v_2)$ in $F$ (if exists). Let $M_2'$ be a perfect matching of $C-v_1$. Then, $M_1'\cup M_2'\cup\{v_1v_2,xy\}$ is a rainbow matching of size $k+1$, a contradiction. This completes the proof. \qed

\section{Proof of Theorem \ref{main1}}

\noindent \textbf{Proof of Theorem \ref{main1}.} By equation (\ref{low-upp}), we have $ar(n,F_{k+1})\geq \ex(n,F_{k})+2$. Therefore the lower bound follows. For $k=1$, by Theorem \ref{Mantel},
we can derive that
$$
\lfloor\dfrac{n^2}{4} \rfloor+2=\ex(n,K_3)+2 \leq ar(n,F_2) \le \ex(n,F_2)+1=(\lfloor \dfrac{n^2}{4} \rfloor+1)+1.
$$
Thus, $ar(n,F_2)= \lfloor\frac{n^2}{4} \rfloor+2$ for $n \geq 5$. Therefore, we can assume $k \geq 2$. Let $c(n,k):=\ex(n,F_k)+2$.


For the upper bound, we prove it by contradiction. Suppose that the result does not hold. Then there exists an edge-coloring $c:E(K_n)\rightarrow [c(n,k)]$ such that $c$ is surjective and the colored graph denoted by $H$ contains no rainbow $ F_{k+1} $. 
Let $G$  be a rainbow subgraph with exactly $c(n,k)$ colors.  Clearly, $G$ is $F_{k+1}$-free and $\Delta(G)>n/2$. We discuss two cases.

\medskip
\textbf{Case 1.}~$\delta(G) > (n/2)-(k+1)$.
\medskip

By Lemma \ref{Gpartition}, there exists a partition $ V(G) =V_0\cup V_1$ such  that $V_0\neq \emptyset$, $V_1\neq \emptyset$ and equations (\ref{5}) and (\ref{6}) hold. By Lemma \ref{y}, we  can derive that
\begin{align*}
|V_0|\cdot|V_1|&\geq e(G)-f(k,k)\\
&= \ex(n,F_k)+2-f(k,k)\\
&=\begin{cases}
	\lfloor \dfrac{n^2}{4} \rfloor-\dfrac{3}{2}k+\dfrac{5}{2},& \mbox{if $k$ is odd},\medskip\\
	\lfloor \dfrac{n^2}{4} \rfloor-\dfrac{5}{2}k+2,& \mbox{if $k$ is even}.
\end{cases}
\end{align*}
i.e.,
\begin{align}\label{V0V1-low}
|V_0|\cdot|V_1|&\geq \begin{cases}
	\lfloor \dfrac{n^2}{4} \rfloor-\dfrac{3}{2}k+\dfrac{5}{2},& \mbox{if $k$ is odd},\medskip\\
	\lfloor \dfrac{n^2}{4} \rfloor-\dfrac{5}{2}k+2,& \mbox{if $k$ is even}.
\end{cases}
\end{align}
Let $|V_0|=\lceil n/2\rceil+r$ and $|V_1|=\lfloor n/2\rfloor-r$, where $r$ is an integer. By equation (\ref{V0V1-low}), one can see that
\begin{align}\label{r-upp}
r^2&\leq \begin{cases}
	\dfrac{3}{2}k-\dfrac{5}{2},& \mbox{if $k$ is odd},\medskip\\
	\dfrac{5}{2}k-2,& \mbox{if $k$ is even}.
\end{cases}
\end{align}
Let $G_0=G[V_0]$ and $G_1=G[V_1]$.  Let $G_2$ denote the bipartite graph with bipartition $(V_0,V_1)$ and edge set $E_G(V_0,V_1)$.
By Lemma \ref{Gpartition},
we  have $e(G_i)\leq f(k,k)$ for $i\in \{0,1\}$.  Therefore, we can infer that
\begin{align}\label{eq:16}
e(G_2)\geq e(G)-2f(k,k)\geq  \lfloor n^2/4\rfloor -k^2-7k/2+2.
\end{align}
So for $i\in \{0,1\}$, by equation (\ref{eq:16}) we have the following statements:
\begin{itemize}
\item [(a)] $|V_i'|\geq n/2-k^2-7k/2+2$, where  $V_{i}':=\{v\in V_i\ |\ V_{1-i}=N_{G_2}(v)\}$;
\item [(b)] for any $v\in V_i$, $d_{G_2}(v)\geq n/2-k^2-7k/2+2$.
\end{itemize}

For $i\in \{0,1\}$, define $S_i=\{ x \in V_i'\ |\ d_{G_i}(x)=0\}$ and $s_i=|S_i|$. That is, $S_i$ is a set of all vertices in $V_i'$ which are not incident with every edge in $G_i$. One can see that
\begin{align*}
\ |S_i|  &\geq |V_i'| - 2(k^2-k+1) \\&\geq \frac{n}{2}-3k^2-\frac{3}{2}k>8k^2.
\end{align*}
For $i\in \{0,1\}$, since $H[S_i]$ contains no rainbow matchings of size $k+1$ and
\begin{align}\label{same-neigh}
\ \frac{|S_i|}{2k}  &\geq  \frac{\frac{n}{2}-3k^2-\frac{3}{2}k}{2k}>4k,
\end{align}
$H[S_i]$ has a monochromatic matching $\mathcal{M}_i$ on color $c_i$ with at least $4k$ edges. Let $G'$ be the rainbow subgraph obtained by removing the  edges  colored by $c_0$ and $c_1$ from $G$. Note that
\begin{align}\label{e(G')-low}
e(G')\geq ex(n,F_k).
\end{align}

If $H[V_i]$ contains a rainbow $K_{1,k+1}$ or a rainbow $(k+1)K_2$, then by  (a) and (b), both of $K_{1,k+1}$ and  rainbow $(k+1)K_2$ can be greedily extended into a rainbow $F_{k+1}$ in $H$, which is a contradiction. Therefore,   we can assume that
for $i\in \{0,1\}$, $H[V_i]$ contains  neither rainbow $K_{1,k+1}$ nor rainbow $(k+1)K_2$.
Furthermore, we can  assume that $\nu(G_0)\geq 1$ and $\nu(G_1)\geq 1$. If not, suppose $\nu(G_i)=0$.  Let $u_1,u_2\in V(G_i)$ and let $v_1v_2\in E(G)$ such that $c(u_1u_2)=c(v_1v_2)$. Then $F=G-v_1v_2+u_1u_2$ is a desired rainbow subgraph of $H$ with $c(n,k)$ edges.

We claim that $G'$ contains no $F_k$. On the contrary, suppose that $G'$ contains an $F_k$ with center vertex $u$ in $V_i$, say $V_0$. Since $\mathcal{M}_1$ is a monochromatic matching of size at
least $4k$ and $V(\mathcal{M}_1)\subseteq S_1$, by the definition of $S_1,$  there exists an edge $v_1v_2\in \mathcal{M}_1$ such that $uv_1,uv_2\in E(G')$. By the definition of $G'$, $c(v_1v_2)\notin \{c(e)\ |\ e\in E(G')\}$, meaning there is no edge in $G'$ with color $c(v_1v_2)$. Thus,
$E(F_k)\cup \{v_1v_2,uv_1,uv_2\}$ induces a rainbow $F_{k+1}$ of $H$ with center $u$, which is a contradiction.

By Theorem \ref{Fkextremal} and equation (\ref{e(G')-low}), we can deduce that $e(G') = ex(n,F_k)$,
which implies that $c_0\neq c_1$, that is, $c_0=c(f_0)\neq c(f_1)=c_1$, for any $f_0\in \mathcal{M}_0$ and $f_1\in \mathcal{M}_1$.  Since $G' \in EX(n,F_k)$, by Lemma \ref{extremalgraph},  either $e(G'[V_0])=ex(n,\{K_{1,k},kK_2\})$ or $e(G'[V_1])=ex(n,\{K_{1,k},kK_2\})$ must hold. Without loss of generality, assume that $e(G'[V_0])=ex(n,\{K_{1,k},kK_2\})$.  
 Note that  $G'[V_0]$ contains a matching denoted by $M_0$ of size $k-1$. Recall that $|\mathcal{M}_i|>4k$ for $i\in \{0, 1\}$. Therefore, we can select $e_0\in \mathcal{M}_0$ such that $M_0\cup \{e_0\}$ is a rainbow matching of $H[V_0]$. Since $V(\mathcal{M}_1)\subseteq S_1$, by the definition of $S_1$,   we can choose $x_1x_2\in \mathcal{M}_1$ such that for $i\in \{1,2\}$, $x_i$ is adjacent to all vertices of $V_0$ in $G'$. Let $v\in S_0-V(M_0\cup \{e_0\})$. Then $M_0\cup \{e_0,x_1x_2,vx_1,vx_2\}\cup \{x_1y\ |\ y\in V(M_0\cup \{e_0\})\}$ induces a  rainbow $F_{k+1}$ of $H$ with center $x_1$, which is a contradiction.

\medskip
\textbf{Case 2.}~$\delta(G) \le (n/2)-(k+1)$.
\medskip

Define $\mathcal{G}_0:=G$. Let $x_0\in V(\mathcal{G}_0)$ such that $d_{\mathcal{G}_0} (x_0) =\delta(\mathcal{G}_0) \le (n/2)-(k+1)$. Define $\mathcal{G}_1=\mathcal{G}_0 - x_0$. Note  the following calculation for the number of edges in $\mathcal{G}_1$:
\begin{align*}
e(\mathcal{G}_1)  &= e(\mathcal{G}_0)-\delta(\mathcal{G}_0) \\
& > (\frac{n^2}{4}-1) -\frac{n}{2} +(k+1) +(k^2-\frac{3}{2}k+2) \\
& =\frac{(n-1)^2}{4}+(k+\frac{3}{4})+(k^2-\frac{3}{2}k+2)-1.
\end{align*}
If there exists a vertex $x_{1} \in V(\mathcal{G}_1)$ such that $d_{\mathcal{G}_1} \left( x_{1} \right) \le (n-1)/2-(k+1)$, then define  $\mathcal{G}_2=\mathcal{G}_1- x_{1}$. Continue this process as long as $\delta(\mathcal{G}_t) \le (|\mathcal{G}_t|/2)-(k+1)$. After $r$ steps,   we obtain a subgraph $\mathcal{G}_r$ with $\delta(\mathcal{G}_r) > (|\mathcal{G}_r|/2)-(k+1)$. By induction, one can show that
\[
e(\mathcal{G}_r) > \frac{(n-r)^2}{4}+r(k+\frac{3}{4})+(k^2-\frac{3}{2}k+2)-1.
\]
On the other hand, since $G$ is $F_{k+1}$-free, $\mathcal{G}_{r}$ is also $F_{k+1}$-free. By Lemma \ref{aa}, we can assume $e(\mathcal{G}_{r}) \le \frac{(n-r)^2}{4}+k(n-r)$. Hence, we deduce that $n>2r$ and thus infer that $n-r> n/2 \ge 24(k+1)^2$. Following the same discussion as Case 1,
$H[V(\mathcal{G}_r)]$ contains a rainbow $F_{k+1}$, which is a contradiction. 		
This completes the proof  of Theorem \ref{main1}. \qed

\end{document}